\documentclass[11pt,leqno]{article} 
\usepackage{graphics}
\newtheorem{thm}{Theorem}[section]
\newtheorem{lma}{Lemma}[section]
\newtheorem{cor}{Corollary}[section]
\newtheorem{prop}{Proposition}[section]
\newtheorem{rmk}{Remark}[section]
\newcommand{\beqa}{\begin{eqnarray}}
\newcommand{\eeqa}{\end{eqnarray}}

\newcommand{\pf}{\noindent {\bf Proof:} $\s$ }
\newcommand{\epf}{ \hfill$\diamondsuit$ \medskip}

\newcommand{\beq}{\begin{equation}}
\newcommand{\eeq}{\end{equation}}
\newcommand{\lbl}{\label}
\newcommand{\s}{\; \;}

\newcommand{\ep}{\epsilon}

\newcommand{\la}{\lambda}

\newcommand{\ra}{\rightarrow}
\newcommand{\al}{\alpha}

\newcommand{\p}{\varphi}
\title{Global solution curves in harmonic parameters, and multiplicity of solutions}

\author{
Philip Korman   \\ 
Department of Mathematical Sciences \\ 
University of Cincinnati \\ 
Cincinnati Ohio 45221-0025 \\
}

\date{}

\begin{document}

\maketitle
\begin{abstract} 
For  the semilinear Dirichlet problem 
\[
\Delta u+g(u)=f(x)  \s \mbox{for $x \in \Omega$}, \s u=0 \s \mbox{on $\partial \Omega$}
\]
decompose $f(x)=\mu _1 \p _1+e(x)$, where $\p _1$ is the principal  eigenfunction of the Laplacian with zero boundary conditions, and $e(x) \perp \p _1$ in $L^2(\Omega)$, and similarly write $u(x)=   \xi _1 \p _i+U (x)$, with $ U \perp \p _1$  in $L^2(\Omega)$. We study properties of the solution curve $(u(x),\mu _1)(\xi _1)$, and in particular its section $\mu _1=\mu _1(\xi _1)$, which governs the multiplicity of solutions. We consider both general nonlinearities, and some important classes of equations, and obtain detailed description of solution curves under the assumption $g'(u)<\la _2$. We obtain particularly detailed results in case of one dimension.
This approach is well suited for numerical computations, which we perform to illustrate our results.
 \end{abstract}

\begin{flushleft}
\end{flushleft}

\begin{flushleft}
\end{flushleft}

\section{Introduction}
\setcounter{equation}{0}
\setcounter{thm}{0}
\setcounter{lma}{0}
\setcounter{rmk}{0}
\setcounter{prop}{0}
\setcounter{cor}{0}

We consider the semilinear Dirichlet problem
\beq
\lbl{intro1}
\s\s \Delta u+g(u)=f(x)  \s \mbox{for $x \in \Omega$}, \s u=0 \s \mbox{on $\partial \Omega$} \,,
\eeq
 on a smooth and bounded domain $\Omega \in R^m$. We are interested in the existence, multiplicity and numerical computation of solutions. 
On the surface the problem (\ref{intro1}) does not contain any parameters. It turns out there are parameters (somewhat hidden), which allow one to understand the solution set. These are the harmonics, particularly the first harmonic, of $f(x)$ and of the solution $u(x)$.
\medskip

Represent the given function $f(x)\in L^2(\Omega)$ by its Fourier series $f(x)=\mu _1 \p _1+\cdots +\mu _n \p _n+e(x)$
where 
$\p _k$ is the $k$-th eigenfunction of the Laplacian on $\Omega$ with zero boundary conditions, and $e(x) \perp \p _k$ in $L^2(\Omega)$, $k=1, \ldots,n$. Likewise decompose  the solution in the form $u(x)= \Sigma _{i=1}^n  \xi _i \p _i+U (x)$, with $ U \perp \p _k$, $k=1, \ldots,n$  in $L^2(\Omega)$. In \cite{K1} and \cite{K9} we studied the map $\xi=\left(\xi _1,\xi _2, \ldots , \xi _n \right) \ra \mu=\left(\mu _1,\mu _2, \ldots , \mu _n \right)$ which determines the multiplicity of solutions of (\ref{intro1}). Indeed, if at some $\mu _0$ this map has $k$ pre-images, then  the problem (\ref{intro1}) has exactly $k$ solutions at $\mu= \mu _0$. In this paper we concentrate on the {\em solution curves}, when the parameters $\xi _i$ are varied one at a time. Particularly, we shall study the curve $\mu _1=\mu _1(\xi _1)$, which we call the {\em principal solution curve}. We shall also consider the curves $\mu _k=\mu _k(\xi _k)$ for a class of oscillatory problems.
\medskip

We now describe one of our results which demonstrates some of the advantages of working with curves of solutions. For a class of  problems
\beq
\lbl{intro2}
\s\s \Delta u+g(u)=\mu _1 \p _1+e(x) \s \mbox{for $x \in \Omega$}, \s u=0 \s \mbox{on $\partial \Omega$} \,,
\eeq
with the solution decomposed as $u(x)=  \xi _1 \p _1+U (x)$ ($U(x) \perp \p _1$), we showed the  existence of  a continuous solution curve $(u(x),\mu _1)(\xi _1)$  which {\em exhausts the solution set of  (\ref{intro2})}. A section of this curve $\mu _1=\mu _1(\xi _1)$ (which is a planar curve) will determine the multiplicity of solutions. (If the value of $\mu _1^0$ is achieved at $k$ values of $\xi _1$, the problem  (\ref{intro2}) has $k$ solutions at $\mu _1=\mu _1^0$). Assuming that $\frac{g(u)}{u}$ crosses the first eigenvalue on the interval $(-\infty,\infty)$ (see the Theorem \ref{thm:62} below) it was possible to show that $\mu _1(\xi _1) \ra \infty$ as $\xi _1 \ra \pm \infty$. It follows that the continuous function $\mu _1(\xi _1)$ has a global minimum value $\mu _0$, and it takes on  at least twice any value in $(\mu _0,\infty)$. We conclude that the problem (\ref{intro2}) has either zero, at least one, or at least two solutions depending on whether $\mu _1< \mu _0$,  $\mu _1= \mu _0$ or  $\mu _1> \mu _0$ respectively, which is  an extension of the well known result of H. Amann and P. Hess \cite{AH}, since we do not require the limits $\lim _{u \ra \pm \infty} \frac{g(u)}{u}$ to exist. The classical result of A. Ambrosetti and G. Prodi \cite{AP}, in the form of M.S. Berger and E. Podolak \cite{BP}, also follows along these lines. We derive similar results if $\frac{g(u)}{u}$ crosses the first eigenvalue on the interval $(0,\infty)$.
A major advantage of our approach is in establishing that all solutions of the problem (\ref{intro2}) lie on a single solution curve, so that they can be computed by a curve following algorithm. Such algorithms (based on Newton's method) are very efficient (fast and accurate), and relatively easy to implement. We present some of our numerical computations.
\medskip

We obtained detailed results for two specific classes of equations.  The first one involves the following model problem
\[
 \Delta u+\la u-u^3=\mu _1 \p _1+e(x) \,, \, \s \mbox{for $x \in \Omega$} \,, \s 
\s u=0 \s \mbox{on $\partial \Omega$} \,,
\]
with a parameter $\la \in (0,\la _2)$, and $\int _{\Omega} e(x) \p _1(x) \, dx=0$. Decompose  $u(x)=\xi _1 \p _1(x)+U(x)$, with $\int _{\Omega} U(x) \p _1(x) \, dx=0$, as above. We show that a typical solution curve $\mu _1=\mu _1(\xi _1)$ is either monotone or $S$-shaped. This problem has been studied previously, see M.S. Berger et al \cite{B} and P.T. Church et al \cite{D} and the references therein, but mostly for the space dimensions $n \leq 4$.
\medskip

The second class involves resonant  problems
\beq
\lbl{july5}
\Delta u +\la _1 u+g(u)=\mu _1 \p _1+e(x) \s \mbox{on $\Omega$}, \s u=0 \s \mbox{on $\partial \Omega$} \,,
\eeq
with $e(x) \in \p _1 ^\perp$ in $L^2(\Omega)$. We wish to find a solution pair $(u, \mu _1)$. For 
{\em bounded} $g(u)$, satisfying $ug(u) \geq 0$ for all $u \in R$, and  $\mu _1=0$,
 D.G. de Figueiredo and W.-M.  Ni \cite{FN} have proved the existence of solutions. R. Iannacci, M.N. Nkashama and J.R. Ward \cite{I} generalized this result to a class of unbounded $g(u)$. In  \cite{K9} we extended the last result to the case $\mu _1 \ne 0$.
Compared with our paper \cite{K9}, here  we do not try to obtain a bound on $\int _{\Omega} |\nabla U|^2 \, dx$ uniformly in $\xi _1$ in general, but either derive such bounds for specific cases, or work around the issue. In particular, that allowed us to drop the technical condition (3.2) of that paper for the problem (\ref{july5}).
\medskip

The most  detailed results are obtained for the one-dimensional case, $n=1$. On the interval $(0,L)$ consider the problem
\beq
\lbl{o1in}
 u''+\la _1 \, u+h(u) \sin u=\mu _1 \p _1(x)+e(x) \,, \s    u(0)=u(L)=0 \,. 
\eeq
Here $\la _1=\frac{\pi ^2}{L^2}$, the principal eigenvalue of $u''$ on $ (0,L)$ corresponding to $\p _1(x)=\sin \frac{\pi}{L} x$, $\mu _1 \in R$, $e(x) \in C (0,L)$ satisfies $\int _0^L e(x) \sin \frac{\pi}{L} x \, dx=0$. Assume that
\beq
\lbl{o2--}
\lim _{u \ra \infty} \frac{h(u)}{u^p}=h_0 \,, \s \mbox{with constants $p \in (0,1)$ and $h_0>0$} \,.
\eeq
It follows from R. Schaaf  and K. Schmitt \cite{SS} that the problem (\ref{o1in}) has infinitely many solutions when $\mu _1=0$. We show that the same is true {\em for all} $\mu _1 \in R$, provided that $p \in (\frac{1}{2},1)$. Moreover, we obtain an asymptotic formula $\mu _1(\xi _1) \sim \frac{2 \sqrt{2}}{\sqrt{\pi \xi _1}} \sin \left( \xi _1-\frac{\pi}{4} \right) h(\xi _1)$ for large $\xi _1$, where $\xi _1$ is the first harmonic of the solution $u(x)$. Our numerical computations suggest that this formula is very accurate, and that it tends to be accurate  for small $\xi _1$ as well (see Figure \ref{fig:1} in Section $6$). In case $p \in (0,\frac{1}{2})$ we showed that
\beq
\lbl{july7}
u(x) \sim \xi _1 \sin \frac{\pi}{L} x+E(x) \,, \s \mbox{for $\xi _1$ large} \,,
\eeq
where $E(x)$ is the unique solution of
\[
u''+\frac{\pi ^2}{L^2} \, u=e(x) \,, \s    u(0)=u(L)=0 \,,  \s \int_0^L u(x) \sin \frac{\pi}{L} x \, dx=0 \,.
\]
The formula (\ref{july7}) gives  {\em an universal asymptotic}, independent of particular $h(u)$. In a forthcoming paper with D.S. Schmidt we establish similar results, and perform similar computations for PDE's.
\medskip

For the resonant problem at higher eigenvalues on $(0,L)$
\[
u''+\frac{k^2 \pi ^2}{L^2} \, u+ \sin u=\mu _k \sin \frac{k\pi}{L} x+e(x) \,,  \s   u(0)=u(L)=0 \,,
\]
where $\int _0^L e(x) \sin \frac{k\pi}{L} x \, dx  =0$,
we decompose the solution in the form  $u(x)=\xi _k \sin \frac{k\pi}{L} x+U(x)$, with  $\int _0^L U(x) \sin \frac{k\pi}{L} x \, dx  =0$,
and prove that all   solutions lie on a unique solution curve $(u(x),\mu _k)(\xi _k)$. Moreover, we  obtain a precise asymptotic formula for  $\mu _k=\mu _k(\xi _k)$, which in particular implies the existence of infinitely many solutions at $\mu _k =0$, see Figure \ref{fig:2}  in Section $6$.

\section{Preliminary results}
\setcounter{equation}{0}
\setcounter{thm}{0}
\setcounter{lma}{0}
\setcounter{rmk}{0}
\setcounter{prop}{0}
\setcounter{cor}{0}

It is well known that on a smooth bounded domain $\Omega \subset R^m$ the eigenvalue problem
\[
\Delta u +\la u=0 \s \mbox{on $\Omega$}, \s u=0 \s \mbox{on $\partial \Omega$}
\]
has an infinite sequence of eigenvalues $0<\la _1<\la _2 \leq \la _3\leq \ldots \ra \infty$, where we repeat each eigenvalue according to its multiplicity, and the corresponding eigenfunctions we denote by $\varphi _k$, and normalize $||\varphi _k||_{L^2(\Omega)}=1$, for all $k$. These eigenfunctions $\varphi _k$ form an orthonormal basis of $L^2(\Omega)$, i.e., any $f(x) \in L^2(\Omega)$ can be written as $f(x)=\Sigma _{k=1}^{\infty} a_k \p _k$, with the series convergent in $L^2(\Omega)$, see e.g., L. Evans \cite{E}.  The following lemma is standard.

\begin{lma}\lbl{lma:1}
Assume that $u(x) \in L^2(\Omega)$, and $u(x)=\sum _{k=n+1}^{\infty} \xi _k \p _k$. Then
\[
\int_\Omega |\nabla u|^2 \, dx  \geq \la _{n+1} \int _\Omega u^2 \, dx.
\]
\end{lma}

In the following linear problem the function $a(x)$ is given, while $\mu _1, \ldots, \mu _n$, and $w(x)$ are unknown.
The following lemmas were proved in \cite{K1} and  \cite{K9}.
\begin{lma}\lbl{lma:3}
Consider the problem
\beqa
\lbl{9}
& \Delta w+a(x)w=\mu_1 \p _1+ \cdots +\mu _n \p _n, \s  \, \mbox{for $x \in \Omega$}, \\ \nonumber
& w=0 \s \mbox{on $\partial \Omega$}, \\ \nonumber
& \int _\Omega w \p _1 \, dx= \cdots = \int _\Omega w \p _n \, dx=0.
\eeqa
Assume that 
\beq
\lbl{10}
a(x) < \la _{n+1}, \,  \s \mbox{for all $x \in \Omega$}.
\eeq
Then the only solution of (\ref{9}) is $\mu _1 = \cdots =\mu _n=0$, and $w(x) \equiv 0$.
\end{lma}

\begin{cor}\lbl{cor:1}
If one considers the problem (\ref{9}) with $\mu _1 = \cdots = \mu _n =0$, then $w(x) \equiv 0$ is the only solution of that problem.
\end{cor}

\begin{cor}\lbl{cor:2}
With $f(x) \in L^2(\Omega)$, consider the problem
\beqa \nonumber
& \Delta w+a(x)w=f(x) \s\s \mbox{for $x \in \Omega$}\,, \\ \nonumber
& w=0 \s\s \mbox{on $\partial \Omega$}, \\ \nonumber
& \int _\Omega w \p _1 \, dx= \cdots = \int _\Omega w \p _n \, dx=0.
\eeqa
Then there is a constant $c$, so that the following a priori estimate holds
\[
||w||_{H^2(\Omega)} \leq c ||f||_{L^2(\Omega)} \,.
\]
\end{cor}

We shall also use a variation of the above lemma, see \cite{K9}.

\begin{lma}\lbl{lma:4}
Consider the problem ($2 \leq i <n$)
\beqa
\lbl{11}
& \Delta w+a(x)w=\mu_i \p _i+\mu_{i+1} \p _{i+1}+ \cdots +\mu _n \p _n \s \mbox{for $x \in \Omega$}, \\ \nonumber
& w=0 \s \mbox{on $\partial \Omega$}, \\ \nonumber
& \int _\Omega w \p _i \, dx= \int _\Omega w \p _{i+1} \, dx=\cdots = \int _\Omega w \p _n \, dx=0.
\eeqa
Assume that 
\beq
\lbl{12}
\la _{i-1}  \leq a(x) \leq  \la _{n+1}, \,  \s \mbox{for all $x \in \Omega$} \,,
\eeq
with at least one of these inequalities being strict.
Then the only solution of (\ref{11}) is $\mu _i = \cdots =\mu _n=0$, and $w(x) \equiv 0$.
\end{lma}

\begin{cor}\lbl{cor:3}
If one considers the problem (\ref{11}) with $\mu _i = \cdots = \mu _n =0$, then $w(x) \equiv 0$ is the only solution of that problem. Consequently, for the problem 
\beqa \nonumber
& \Delta w+a(x)w=f(x) \s\s \mbox{for $x \in \Omega$}\,, \\ \nonumber
& w=0 \s\s \mbox{on $\partial \Omega$}, \\ \nonumber
& \int _\Omega w \p _i \, dx= \cdots = \int _\Omega w \p _n \, dx=0.
\eeqa
there is a constant $c$, so that the following a priori estimate holds
\[
||w||_{H^2(\Omega)} \leq c ||f||_{L^2(\Omega)} \,.
\]
\end{cor}

\section{Continuation of solutions}
\setcounter{equation}{0}
\setcounter{thm}{0}
\setcounter{lma}{0}
\setcounter{rmk}{0}
\setcounter{prop}{0}
\setcounter{cor}{0}

Any $f(x) \in L^2(\Omega)$ can be decomposed as $f(x)=\mu_1 \p _1+ \cdots +\mu _n \p _n+e(x)$, with $e(x)$ orthogonal to $\p _1,  \ldots, \p _n$ in $L^2(\Omega)$. We consider the following  boundary value problem
\beqa
\lbl{2}
& \Delta u+g(u)=f(x)=\mu_1 \p _1+ \cdots +\mu _n \p _n+e(x) \s \mbox{for $x \in \Omega$} \,, \\ \nonumber
& u=0 \s \mbox{on $\partial \Omega$} \,.
\eeqa
If $u(x) \in H^2(\Omega) \cap H^1_0(\Omega)$ is a solution of (\ref{2}), we decompose it likewise  as 
\beq
\lbl{5}
u(x)= \Sigma _{i=1}^n  \xi _i \p _i+U (x),
\eeq
where $U (x)$ is orthogonal to $\p _1, \ldots, \p _n$ in $L^2(\Omega)$. 
We pose the following inverse problem: keeping $e(x)$ fixed, find $\mu=\left( \mu _1, \ldots, \mu _n \right)$ so that the problem (\ref{2}) has a solution for any prescribed  $\xi=\left( \xi _1, \ldots, \xi _n \right)$. Under the conditions given below this problem has a unique solution, and therefore  we shall call $\xi=\left( \xi _1, \ldots, \xi _n \right)$ the {\em $n$-signature of the solution}.
\medskip

The following result generalizes the corresponding one in \cite{K9}, by dropping  the technical condition (3.2) of that paper.

\begin{thm}\lbl{thm:1}
For the problem (\ref{2}) assume that $g(u) \in C^2(R)$, $f(x) \in L^2(\Omega)$, and
\beq
\lbl{3}
g'(u)<\la _{n+1} \,, \s \mbox{for all $u \in R$} \,,
\eeq
\beq
\lbl{3a}
\s\s\s \s |g(u)|<\gamma |u| +c  \,, \s \mbox{with constants $0<\gamma<\la _{n+1}$, $c \geq 0$, and   $u \in R$} \,.
\eeq
Then given any $\xi=\left( \xi _1, \ldots, \xi _n \right)$, one can find a unique $\mu=\left( \mu _1, \ldots, \mu _n \right)$ for which the problem (\ref{2}) has a solution $u(x) \in H^2(\Omega) \cap H^1_0(\Omega)$ of $n$-signature $\xi$. This solution is unique. 
\end{thm}

\pf
We embed the problem (\ref{2}) into a family of problems
\beqa
\lbl{2k}
& \Delta u+kg(u)=\mu_1 \p _1+ \cdots +\mu _n \p _n+e(x) \s \mbox{for $x \in \Omega$}, \\ \nonumber
& u=0 \s \mbox{on $\partial \Omega$} \,,
\eeqa
depending on a parameter $0 \leq k \leq 1$.
Decompose $e(x)=\sum _{j=n+1}^{\infty} e_j \p _j$. When $k=0$, the problem (\ref{2k}) has infinitely many solutions. The unique solution of (\ref{2k}) with signature $\xi$ is $u(x)=\sum _{j=1}^{n} \xi _j \p _j-\Sigma _{j=n+1}^{\infty} \frac{e_j}{\la _j} \p _j$, corresponding to $\mu _j =-\la _j \xi _j$, $j=1, \ldots, n$. We shall use the implicit function theorem to continue this solution in $k$, obtaining a curve $(u(x),\mu)(k)$ (with the $n$-signature of $u(x)$ being fixed at $\xi$). Writing $u(x)= \sum _{i=1}^n  \xi _i \p _i+U (x)$, we multiply the equation (\ref{2k}) by $\p _i$, and integrate
\beq
\lbl{16}
\mu _i=-\la _i \xi _i+ k\int _\Omega g \left( \sum _{i=1}^n  \xi _i \p _i+U \right) \p _i \, dx, \s i=1, \ldots, n \,.
\eeq
Using these expressions in (\ref{2k}), obtain
\beq
\lbl{17}
\s \s \s \Delta U+kg\left( \sum _{i=1}^n  \xi _i \p _i+U \right)-k \sum _{i=1}^n \int _\Omega g \left( \sum _{i=1}^n  \xi _i \p _i+U \right)  \p _i \, dx \, \p _i=e(x) \,, 
\eeq
\[
 U=0 \s \mbox{on $\partial \Omega$} \,.
\]
The equations (\ref{16}) and (\ref{17}) constitute the classical Lyapunov-Schmidt decomposition of our problem (\ref{2}).
Define $H^2_{{\bf 0}}$ to be the subspace of $H^2(\Omega) \, \cap H^1_0(\Omega)$, consisting of functions with zero $n$-signature:
\[
H^2_{{\bf 0}} =  \left\{ u \in H^2(\Omega) \cap H^1_0(\Omega) \; | \; \int _\Omega u \p _i \, dx =0, \; i=1, \ldots, n \right\}.
\]
We recast the problem (\ref{17}) in the operator form as
\[
F(U, k) =e(x),
\]
where $ F(U, k) : H^2_{{\bf 0}}  \times R  \ra L^2(\Omega)$ is given by the left hand side of (\ref{17}). Compute the Frechet derivative
\[
F_{U}(U, k)w=\Delta w+kg' \left( \sum _{i=1}^n  \xi _i \p _i+U \right)w-\mu^*_1 \p _1-\cdots -\mu^* _n \p _n \,,
\]
where we denoted $\mu^*_i=k \int _\Omega g' \left( \sum _{i=1}^n  \xi _i \p _i+U \right) w \p _i \, dx$. By Lemma \ref{lma:3} the map $F_{U}(U, k)$ is injective. Since this map is Fredholm of index zero, it is also surjective. The implicit function theorem applies, giving us locally a curve of solutions $U=U(k)$. Then we compute $\mu=\mu (k)$ from (\ref{16}).
\medskip

To show that this curve  $(u(k),\mu (k))$ can be continued for all $k$,  one needs to show that it cannot go to infinity at some $k \in (0,1)$, i.e., one  needs an a priori estimate. Since the $n$-signature of the solution is fixed, we only need to estimate $U$. We claim that there is a constant $c>0$, $c=c\left( \xi _1, \ldots, \xi _n \right)$, so that 
\beq
\lbl{18}
||U||_{H^2(\Omega)} \leq c \,.
\eeq
Using projections rewrite the equation in (\ref{17}) as
\beq
\lbl{18a}
\Delta U+kPg\left( \sum _{i=1}^n  \xi _i \p _i+U \right) =e(x) \,,
\eeq
where $P$ is the projection on $\{\p _1, \ldots , \p _n \}^{\perp}$ in $L^2$.
Multiply (\ref{18a}) by $U$ and integrate:
\beq
\lbl{19a}
\s\s\s -\int _{\Omega} |\nabla U|^2 \, dx+k \int _{\Omega} Pg\left( \sum _{i=1}^n  \xi _i \p _i+U \right) U  \, dx =\int _{\Omega} U(x) e(x)  \, dx \,.
\eeq
Estimate (with $u=\sum _{i=1}^n  \xi _i \p _i+U$)
\[
|k \int _{\Omega} Pg(u) U  \, dx | \leq ||Pg(u)||_{L^2} ||U||_{L^2} \leq ||g(u)||_{L^2} ||U||_{L^2} \,.
\]
By (\ref{3a}) it follows that
\[
g^2(u)<(\gamma+\ep)^2 |u|^2 +c_1 \,,
\]
with some small $\ep >0$, $0<\gamma<\la _{n+1}$, and $c_1>0$. Then
\[
||g\left( \sum _{i=1}^n  \xi _i \p _i+U \right)||_{L^2}^2 <(\gamma+\ep)^2 ||U||_{L^2}^2+c_2 \,,
\]
with some $c_2>0$, and using Lemma \ref{lma:1} we conclude from (\ref{19a}) an estimate on $||U||_{L^2}$. Then the estimate (\ref{18}) follows from (\ref{18a}), by using (\ref{3a}) and elliptic regularity.
\medskip

Finally, if the problem (\ref{2}) had a different solution $(\bar u,\bar \mu )$ with the same signature $\xi$, then starting with (\ref{2k}) at $k=1$, we would continue the solution  back in $k$, obtaining at $k=0$ a different  solution of the linear problem of signature $\xi$ (since solution curves do not intersect by the implicit function theorem), which is impossible.
\epf

\begin{rmk}
\lbl{rmk:1}
Observe that the condition (\ref{3a}) follows from (\ref{3}), provided either one of the following two conditions holds:
\beq
\lbl{3b}
(u-a)g(u)>0 \,, \s \mbox{for some $a \in R$, and all $u \in R$} \,,
\eeq
or
\beq
\lbl{3c}
\mbox{$g(u)$ is increasing for $|u|$ large} \,.
\eeq
\end{rmk}

The Theorem \ref{thm:1} implies that the value of $\xi =(\xi_1, \ldots, \xi_n)$ uniquely identifies the solution pair $( u(x),\mu)$, where $\mu =(\mu_1, \ldots, \mu _n)$. Hence, the solution set of (\ref{2}) can be completely described by the map: $\xi \in R^n \ra \mu \in R^n$, which we call the {\em solution manifold}.  We show next that the solution manifold is connected.

\begin{thm}\lbl{thm:2}
Under the conditions (\ref{thm:3}),(\ref{3a}) of Theorem \ref{thm:1}, the solution $(u,\mu_1,\dots,\mu _n)$ of (\ref{2}) is a continuous function of $\xi=(\xi_ 1, \dots ,\xi _n)$. Moreover,  we can continue solutions of any signature $\bar \xi$ to solution of arbitrary signature $\hat \xi $ by following any bounded continuous curve in $R^n$ joining $\bar \xi$ and $\hat \xi$. 
\end{thm}

\pf
We use the implicit function theorem to show that any solution of (\ref{2}) can be continued in $\xi$. The proof is essentially the same as for the continuation in $k$ above. After performing the same Lyapunov-Schmidt decomposition,
we recast the problem (\ref{17}) in the operator form
\[
F(U,\xi)=e(x) \,,
\]
where $F \, : \, H^2_{{\bf 0}} \times R^n  \ra L^2$ is defined  by the left hand side of (\ref{17}). The Frechet derivative $F_{U}(U, \xi)w$ is the same as before, and by the implicit function theorem we have locally $U=U(\xi)$. Then we compute $\mu=\mu (\xi)$ from (\ref{16}). We use the same a priori bound (\ref{18}) to continue the curve for all $\xi \in R^n$. (The bound (\ref{18}) is uniform, once the curve  joining $\bar \xi$ and $\hat \xi$ is fixed.)
\epf

If the conditions (\ref{thm:3}),(\ref{3a}) hold with $n=1$, in other words 
\beq
\lbl{pc}
g'(u)<\la _2 \,, \s \mbox{ for all $u$} \,, 
\eeq
and the condition (\ref{2k1*}) below holds, we conclude by the Theorem \ref{thm:1} that the problem
\beq
\lbl{2k1}
\Delta u+g(u)=\mu_1 \p _1+ e(x) \s\s \mbox{for $x \in \Omega$} \,, \s
 u=0 \s \mbox{on $\partial \Omega$} \,,
\eeq 
has a continuous solution curve $(u(x),\mu _1)(\xi _1)$, which we call {\em the principal solution curve}, with $\xi _1 \in R$ serving as a global parameter. We now obtain a solution curve for (\ref{2k1}) without the assumption (\ref{pc}).

\begin{prop}\lbl{prop:n1}
For the problem (\ref{2k1}) assume that
\beq
\lbl{2k1*}
\s\s\s\s |g(u)|<\gamma |u| +c  \,, \s \mbox{with constants $\gamma<\la _2$, $c \geq 0$, and all  $u \in R$} \,.
\eeq
Then for any $\xi _1 \in R$ one can find a  $\mu _1$ for which the problem (\ref{2k1}) has a  solution $u(x)$, with the first harmonic equal to $\xi _1$.
\end{prop}

Observe that again we have a  solution curve $(u(x),\mu _1)(\xi _1)$,  although this time the curve is not necessarily continuous. Also now $\xi _1$ is   not necessarily a global parameter (there could be other  solution curves).
\medskip

\pf
Substitution of $u=\xi _1 \p _1+U$ into  (\ref{2k1}) gives
\beq
\lbl{20m1}
\s\s\s\s -\la _1 \xi _1 \p _1+\Delta U+g(\xi _1 \p _1+U)=\mu_1 \p _1+ e \s\s \mbox{$x \in \Omega$} \,, \s
 U=0 \s \mbox{on $\partial \Omega$} \,.
\eeq 
Multiplying this equation by $\p _1$ and integrating one gets an expression for $\mu _1$, which is substituted  back into (\ref{20m1}) obtaining
\beqa
\lbl{20m2}
& \mu _1=-\la _1 \xi _1+\int _{\Omega} g(\xi _1 \p _1+U) \p _1 \, dx \\ \nonumber
& \Delta U+Pg(\xi _1 \p _1+U)= e \s\s \mbox{$x \in \Omega$} \,, \s
 U=0 \s \mbox{on $\partial \Omega$} \,. \nonumber
\eeqa
Here $Pg(\xi _1 \p _1+U)=g(\xi _1 \p _1+U)-\p _1 \int _{\Omega} g(\xi _1 \p _1+U) \p _1 \, dx$ gives the projection of $g(\xi _1 \p _1+U)$ on the subspace  $\p _1^{\perp}$ in $L^2(\Omega)$.
The equations in (\ref{20m2}) constitute the classical Lyapunov-Schmidt reduction.
Proceeding as in the proof of  Theorem \ref{thm:1} (the argument that begins at (\ref{19a})), we get an a priori bound on $||U||_{L^2}$, using the condition (\ref{2k1*}). Define the space $X=\p _1^{\perp}$ in $L^2(\Omega)$, and the operator $T \, : \, X \ra X$ as $T(U)=\Delta ^{-1} \left(e(x)-Pg(\xi _1 \p _1+U) \right)$. The operator $T$ is continuous and compact, and the set $\{ U \in X \, {\Huge |} \, U=\la T(U) \s \mbox{for some $0 \leq \la \leq 1$} \}$ is bounded. By  Schaefer's fixed point theorem, see e.g., \cite{E}, the second equation in (\ref{20m2}) has a solution $U(x)$. Then $\mu _1$ is determined from the first equation in (\ref{20m2}).
\epf

For sublinear $g(u)$ a similar result was obtained in R. Schaaf  and K. Schmitt \cite{SS1}, based on E.N. Dancer \cite{D1}.
\medskip

We now extend the above results. Given a Fourier series $u(x)=\Sigma _{j=1}^{\infty} \xi _j \p _j$, we call the vector $(\xi _i, \ldots, \xi _n)$ to be the $(i,n)$-{\em signature } of $u(x)$, $2 \leq i<n$.
Using Lemma \ref{lma:4} in place of Lemma \ref{lma:3}, we have the following variation of the Theorems \ref{thm:1} and \ref{thm:2}. We decompose  $u(x)=\sum _{j=i}^n \xi _j \p _j(x)+U(x)$ and $f(x)=\sum _{j=i}^n \mu _j \p _j(x)+e(x)$, with $U(x)$ and $e(x)$ orthogonal to $\p _i(x), \ldots , \p _n(x)$ in $L^2(\Omega)$.

\begin{thm}\lbl{thm:3}
For the problem (\ref{2}) assume that the condition (\ref{3a}) hold, and
\[
\la _{i-1} \leq g'(u) \leq \la _{n+1}, \,  \s \mbox{for all $u \in R$} \,,
\]
with at least one of these inequalities being strict.
 Then given any $\xi=\left( \xi _i, \ldots, \xi _n \right)$, one can find a unique $\mu=\left( \mu _i, \ldots, \mu _n \right)$ for which the problem
\beqa
\lbl{20}
& \Delta u+g(u)=\mu_i \p _i+ \cdots +\mu _n \p _n+e(x), \, \s \mbox{for $x \in \Omega$}, \\ \nonumber
& u=0 \s \mbox{on $\partial \Omega$}
\eeqa
 has a solution $u(x) \in H^2(\Omega) \cap H^1_0(\Omega)$ of the $(i,n)$-signature $\xi$. This solution is unique. Moreover, the solution $(u(x),\mu)(\xi)$ is a continuous function of  $\xi$. In addition, we can continue solutions of any  $(i,n)$-signature $\bar \xi$ to solution of arbitrary  $(i,n)$-signature $\hat \xi $ by following any continuous curve in $R^{n-i+1}$ joining $\bar \xi$ and $\hat \xi$.
\end{thm}

\begin{rmk}
\lbl{rmk:*}
In particular, when $i=n$, and $\la _n$ is a simple eigenvalue corresponding to the eigenfuction $\p _n$, this theorem asserts the existence of the unique continuous solution curve $(u(x),\mu _n)(\xi _n)$, for all $\xi _n \in R$,  for the problem
\beqa \nonumber
& \Delta u+g(u)=\mu _n \p _n+e(x), \, \s \mbox{for $x \in \Omega$} \\ \nonumber
& u=0 \s \mbox{on $\partial \Omega$} \,,
\eeqa
provided that
\[
\la _{n-1} \leq g'(u) \leq \la _{n+1}, \,  \s \mbox{for all $u \in R$} \,,
\]
with at least one of these inequalities being strict.
Here $u(x)=\xi _n \p _n(x)+U(x)$, with $U(x)$ and $e(x)$ orthogonal to $\p _n(x)$ in $L^2(\Omega)$.
\end{rmk}

\begin{rmk}
\lbl{rmk:2}
The condition (\ref{3a}) was used only to obtain the estimate (\ref{18}) on $||U||_{H^2(\Omega)}$. Hence, one can drop the condition (\ref{3a}) in  the Theorems \ref{thm:1},   \ref{thm:2} and \ref{thm:3} if   the estimate (\ref{18}) can be obtained in another way.
\end{rmk}

For example, we have the following proposition.

\begin{prop}
\lbl{prop:1}
Consider the problem
\beqa
\lbl{80}
& \Delta u+\la u-u^3=\mu _1 \p _1+e(x) \,, \, \s \mbox{for $x \in \Omega$} \\ \nonumber
& u=0 \s \mbox{on $\partial \Omega$} \,,
\eeqa
with a parameter $\la \in (0,\la _2)$, and $e(x) \in L^2(\Omega)$ satisfying $\int _{\Omega} e(x) \p _1(x) \, dx=0$. Decompose $u(x)=\xi _1 \p _1(x)+U(x)$, with $\int _{\Omega} U(x) \p _1(x) \, dx=0$. Then for each $\xi _1\in R$ one can find a unique solution pair $(u(x),\mu _1)$, and the solution curve $(u(x),\mu _1)(\xi _1)$ is continuous. 
\end{prop}

\pf
The condition (\ref{pc}) holds here. In view of the Remark \ref{rmk:2}, we only need to derive the estimate (\ref{18}), in order to apply Theorem \ref{thm:1}. Multiply the equation (\ref{80}) by $U$ and integrate
\beq
\lbl{81}
- \int _{\Omega}  |\nabla U|^2 \, dx+\la \int _{\Omega} U^2  \, dx - \int _{\Omega} u^3 U  \, dx = \int _{\Omega} Ue \, dx \,.
\eeq
Estimate
\[
\int _{\Omega} u^3 U  \, dx=\int _{\Omega} \left(\xi _1 \p _1+U \right) ^3 U  \, dx=\int _{\Omega} U^4 \, dx+\cdots \geq c_1 \,,
\]
for some constant $c_1=c_1(\xi _1)$. Then from (\ref{81})
\[
\left(\la _2-\la \right)\int _{\Omega}  |\nabla U|^2 \, dx+c_1 \leq \int _{\Omega}  |\nabla U|^2 \, dx-\la \int _{\Omega} U^2  \, dx + \int _{\Omega} u^3 U  \, dx = -\int _{\Omega} Ue \, dx \,,
\]
from which, using the Poincare inequality, we obtain an estimate $\int _{\Omega}  |\nabla U|^2 \, dx \leq c_2$, and then $\int _{\Omega}  U^2 \, dx \leq c_2$, with $c_2=c_2(\xi _1)$. It follows that 
\beq
\lbl{81a}
\int _{\Omega}  |\nabla u|^2 \, dx \leq c_3 \,, 
\eeq
with $c_3=c_3(\xi _1)$. Now multiply (\ref{80}) by $\Delta u$ and integrate
\[
\int _{\Omega} \left(\Delta u \right)^2 \, dx- \la \int _{\Omega} |\nabla u|^2 \, dx+3 \int _{\Omega}  u^2 |\nabla u|^2 \, dx=\int _{\Omega} \Delta u \, e \, dx
\]
which gives an estimate on $\int _{\Omega} \left(\Delta u \right)^2 \, dx$, in view of (\ref{81a}). Then 
we get an estimate on $\int _{\Omega} \left(\Delta U \right)^2 \, dx$, and  using elliptic regularity conclude the estimate (\ref{18}), completing the proof.
\epf

\begin{prop}\label{prop:inf}
Consider the problem
\[
\Delta u +\la u+g(u)=\mu _1 \p _1+e(x) \s \mbox{on $\Omega$}, \s u=0 \s \mbox{on $\partial \Omega$} \,.
\]
Assume that $0 \leq \la <\la _2$,  $\lim _{|u| \ra \infty} \frac{g(uz)}{u}=0$ uniformly in $z \in R$, and $e(x) \perp \p _1$ in $L^2(\Omega)$. Then as $\xi _1 \ra \pm \infty$, we have $\frac{u(x)}{\xi _1} \ra \p _1(x)$ in $H^1(\Omega)$.
\end{prop}

\pf
By the Proposition \ref{prop:n1} we have a solution curve $(u(x),\mu _1)(\xi _1)$, and $\frac{u(x)}{\xi _1}=\p _1(x)+\frac{U(x)}{\xi _1}$. Letting $U=\xi _1 V$ in (\ref{20m1}), obtain
\[
(\la-\la _1)  \p _1+\Delta V+\la V=-\frac{g(\xi _1 \left(\p _1+ V \right))}{\xi _1}+\frac{\mu_1}{\xi _1} \p _1+ \frac{e}{\xi _1}=\frac{\mu_1}{\xi _1} \p _1+o(1) \,.
\]
Multiplying by $V$ and integrating, we conclude that $\int _{\Omega}  |\nabla V|^2 \, dx=o(1)$, as $\xi _1 \ra \pm \infty$.
\epf

\begin{prop}\label{prop:infa}
Consider the problem (with integer $k \geq 2$)
\beq
\lbl{n1-1}
\Delta u +\la _k u+g(u)=\mu _k \p _k+e(x) \s \mbox{on $\Omega$}, \s u=0 \s \mbox{on $\partial \Omega$} \,.
\eeq
Assume that $\la _k$ is a simple eigenvalue, and $\la _{k-1}<\la _k +g'(u)<\la _{k+1}$ for all $u \in R$,  $g(u)$ is bounded  uniformly in $u \in R$, and $e(x) \perp \p _k$ in $L^2(\Omega)$. Then as $\xi _k \ra \pm \infty$, we have $\frac{u(x)}{\xi _k} \ra \p _k(x)$ in $H^2(\Omega)$.
\end{prop}

\pf
By the Remark \ref{rmk:*} we have a solution curve $(u(x),\mu _k)(\xi _k)$, and $\frac{u(x)}{\xi _k}=\p _k(x)+\frac{U(x)}{\xi _k}$.  Multiplying (\ref{n1-1}) by $\p _k$ and integrating over $\Omega$, we conclude a bound on $|\mu _k|$. Letting $u=\xi _k \p _k +U$, followed by $U=\xi _k V$ in (\ref{n1-1}), write the result as 
\[
\Delta V+\la _k V=-\frac{g(\xi _k\left(\p _k+ V \right))}{\xi _k}+\frac{\mu_k}{\xi _k} \p _k+ \frac{e}{\xi _k}=o(1) \,,
\]
as $\xi _k \ra \pm \infty$.
The proof follows by Corollary \ref{cor:3}.
\epf

\section{The principal global solution curve}
\setcounter{equation}{0}
\setcounter{thm}{0}
\setcounter{lma}{0}
\setcounter{rmk}{0}
\setcounter{prop}{0}
\setcounter{cor}{0}

In this section we study the shape of the {\em principal solution curve} $(u(x),\mu _1)(\xi _1)$ of the problem
\beq
\lbl{n1}
\Delta u +g(u)=f(x)=\mu _1 \p _1+e(x) \s \mbox{on $\Omega$}, \s u=0 \s \mbox{on $\partial \Omega$} \,,
\eeq
assuming that for all $u \in R$
\beq
\lbl{n2}
g'(u) \leq \gamma <\la _2   \,,
\eeq
\beq
\lbl{n2a}
|g(u)|<\gamma |u|+c \,, \s \mbox{with constants $0<\gamma <\la _2$, $c \geq 0$} \,.
\eeq
Here $u(x)=\xi _1 \p _1(x)+U(x)$, with $\int _{\Omega} U(x) \p _1(x) \, dx=0$. The existence of this solution curve follows by the Theorem \ref{thm:1}. 
\medskip

We shall use the following results of A.C.  Lazer and  P.J. McKenna \cite{L1} for the problem (\ref{n1}).

\begin{prop}(\cite{L1})\lbl{prop:g1}
Assume that the function $g(u) \,: R \ra R$, $g(u) \in C^1(R)$ satisfies 
\beq
\lbl{g2}
g(u)-\la _1 u \geq c_0|u|-b \,,
\eeq
with constants $c_0>0$ and $b \geq 0$, and $g'(u)$ is bounded on $[0,\infty)$.  Let $f(x) \in C^{\al }(\bar \Omega)$. Denote $M=\max _{\bar \Omega} |f(x)|$.
Then any solution of (\ref{n1}) satisfies $||u||_{C^{1+\al }(\bar \Omega)} \leq c$, with $c=c(M)$.
\end{prop}

\begin{prop}(\cite{L1})\lbl{prop:g2}
Assume that (\ref{n2}) holds. Then any two distinct solution $v$ and $w$ of (\ref{n1}) satisfy $v(x)-w(x) \ne 0$ for all $x \in \Omega$.
\end{prop}

\begin{thm}(\cite{L1})\lbl{thm:g3}
Assume that   (\ref{n2}) holds, and either 
\[
\mbox{(a) \hspace{.4in} $g'(u)$ is strictly increasing in $u$}, or
\]
\[
\mbox{(b) \hspace{.4in} $g'(u)$ is strictly decreasing in $u$} \,.
\]
Then (\ref{n1}) has at most two solutions.
\end{thm}

We have the following general result.

\begin{thm}\lbl{thm:62}
Assume that  the conditions (\ref{n2}) and (\ref{n2a}) hold, and there exist constants $0<\gamma _1 <\la _1<\gamma _2$, and $N>0$ so that 
\beqa
\lbl{g4} 
& \frac{g(u)}{u} <\gamma _1 \,, \s \mbox{for $u<-N$} \,, \\ \nonumber
& \frac{g(u)}{u} >\gamma _2 \,, \s \mbox{for $u>N$} \,. \\ \nonumber
\eeqa
Assume also that $e(x) \in C^{\al }(\bar \Omega)$.
Then all solutions of (\ref{n1}) lie on a unique continuous solution curve $(u(x),\mu _1)(\xi _1)$, and $\mu _1(\xi _1) \ra +\infty$ as $\xi _1 \ra \pm \infty$. Consequently, there exists a constant $\mu _0$ so that the problem (\ref{n1}) has either zero, at least one, and at least two solutions depending on whether $\mu _1< \mu _0$,  $\mu _1= \mu _0$ or  $\mu _1> \mu _0$ respectively. In case $g'(u)$ is strictly increasing for all $u \in R$,  the problem (\ref{n1}) has no solutions for $\mu _1< \mu _0$,  exactly one solution for $\mu _1= \mu _0$, and exactly two solutions for $\mu _1> \mu _0$.
\end{thm}

\pf
The existence and uniqueness of the solution curve $(u(x),\mu _1)(\xi _1)$, or $\mu _1= \mu _1(\xi _1)$ follows by Theorem \ref{thm:1}, we now discuss the properties of this curve. We claim that there exists $\mu ^* \in R$ such that $\mu _1(\xi _1)>\mu ^*$ for all $\xi _1 \in R$. The proof is standard, see A.  Ambrosetti and D. Arcoya \cite{AA}, but we include it for completeness. The condition (\ref{g4}) implies that for some $c \geq 0$ and all $u \in R$ the following two inequalities hold:
\beq
\lbl{g6}
g(u) \geq \gamma _1 u-c \,, 
\eeq
\beq
\lbl{g7}
g(u) \geq \gamma _2 u-c  \,.
\eeq
From (\ref{n1})
\[
\mu _1=\int _{\Omega} g(u) \p _1 \, dx-\la _1 \xi _1 \,.
\]
In case $\xi _1 \geq 0$, we use (\ref{g7}) to get
\[
\mu _1 \geq \gamma _2 \int _{\Omega} u \p _1 \, dx-\la _1 \xi _1-c=\left(\gamma _2-\la _1\right) \xi _1 -c \geq -c.
\]
In case $\xi _1 < 0$, we use (\ref{g6}) to get
\[
\mu _1 \geq \gamma _1 \int _{\Omega} u \p _1 \, dx-\la _1 \xi _1-c=\left(\gamma _1-\la _1\right) \xi _1 -c \geq -c.
\]

Observe that the conditions (\ref{g4}) imply that the condition (\ref{g2}) holds. By Proposition \ref{prop:g1}, $|\mu _1|$ tends to infinity as $\xi _1 \ra \pm \infty$  (if $\mu _1$ is bounded, so is $u(x)$ and hence $\xi _1$ is bounded, a contradiction). Since $\mu _1(\xi _1)$ is bounded from below, $\mu _1(\xi _1) \ra +\infty$ as $\xi _1 \ra \pm \infty$. Let $\mu _0$ be the global minimum value of $\mu _1(\xi _1)$. Then the multiplicity of solutions follows, see Figure \ref{fig:3}. In case $g'(u)$ is strictly increasing, the exact multiplicity count of solutions  follows by Theorem \ref{thm:g3}.
\epf

This theorem provides an extension of the well known result of H. Amann and P. Hess \cite{AH}, since we do not require the limits $\lim _{u \ra \pm \infty} \frac{g(u)}{u}$ to exist (on the other hand, the condition $g'(u)<\la _2$ was not required in \cite{AH}). In case $g'(u)$ is strictly increasing, this theorem recovers the classical result of A. Ambrosetti and G.  Prodi \cite{AP}, in the form of M.S. Berger and E.  Podolak \cite{BP}.
\medskip

\noindent
{\bf Example} We computed the solution curve $\mu _1= \mu _1(\xi _1)$ for the following example (here $u(x)=\xi _1 \sin \pi x+U(x)$)
\beqa 
\lbl{o2aaa}
& \s\s u''+g(u)=\mu _1 \sin \pi x+e(x) \,, \s x \in (0,1) \,, \\ \nonumber
& u(0)=u(1)=0 \,, \nonumber
\eeqa
with $g(u)=\cos u + u \left(\pi ^2+\frac{2}{\pi} \tan ^{-1} u+0.9 \sin \left(\ln (u^2+1) \right) \right)$, $e(x)=\sin 2 \pi x-2\sin 5 \pi x$. Observe that $\int _0^1 e(x)  \sin \pi x \, dx=0$. Here $\la _1=\pi ^2$, $\p _1(x)=\sin \pi x $, $\la _2=4\pi ^2$, and one checks that the Theorem \ref{thm:62} applies, while the above mentioned result of H. Amann and P. Hess \cite{AH} does not (the limits $\lim _{u \ra \pm \infty} \frac{g(u)}{u}$ do not exist). The solution curve $\mu _1= \mu _1(\xi _1)$  is presented in Figure \ref{fig:3}. 
Here, and in Figures \ref{fig:1} and \ref{fig:2}, we used a program written jointly with D.S. Schmidt, see \cite{KS} for a detailed explanation of this program.

\begin{figure}
\begin{center}
\scalebox{0.70}{\includegraphics{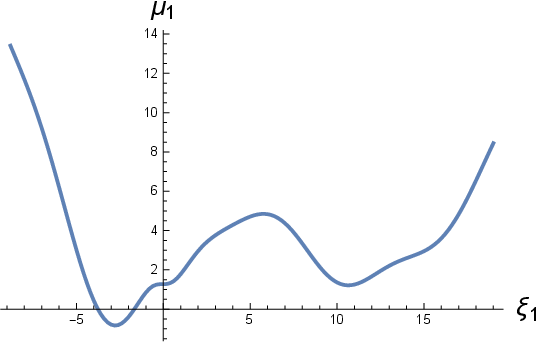}}
\end{center}
\caption{ The  solution curve $\mu _1= \mu _1(\xi _1)$ of (\ref{o2aaa})}
\lbl{fig:3}
\end{figure}
\medskip

We shall  need the following extension of a lemma from \cite{Knew}. 
\begin{lma}\lbl{lma:6}
Let $u(x)$ be a solution of the problem
\beq
\lbl{22}
\Delta u +\la _1 u+a(x)u=\mu _1 \p _1+e(x) \s \mbox{on $\Omega$}, \s u=0 \s \mbox{on $\partial \Omega$},
\eeq
with $e(x) \in \p _1 ^\perp$ in $L^2(\Omega)$, and $a(x) \in C(\Omega)$. Assume  there is a constant $\gamma$, so that
\[
0 \leq a(x)\leq \gamma <\la_2- \la _1, \s\s \mbox{for all $x \in \Omega$} \,.
\]
Write the solution of (\ref{22}) in the form $u(x)=\xi _1 \p _1+U$, with $U \in \p _1 ^\perp$, and assume that
\beq
\lbl{22a}
\xi _1 \mu _1 \leq 0 \,.
\eeq
Then there exists a constant $c_0$, so that
\beq
\lbl{22.2}
\int_\Omega |\nabla U|^2 \, dx \leq c_0 \,, \s\s \mbox{uniformly in $\xi _1 $ satisfying (\ref{22a})}\,.
\eeq
If, in addition, one has 
\beq
\lbl{22.2a}
a(x)>\ep >0 \, \s\s \mbox{for all $x \in \Omega$} \,,
\eeq
for some constant $\ep$, then the condition (\ref{22a}) may be replaced by
\beq
\lbl{22.2b}
\ep \xi _1^2-\xi _1 \mu _1>0 \,.
\eeq
\end{lma}

\pf
Substitution of $u=\xi _1 \p _1+U$ into (\ref{22}) gives
\beq
\lbl{22.1}
\s \s\s\s \Delta U +\la _1 U+a(x)\left(\xi _1 \p _1+U \right)=\mu _1 \p _1+e(x) \s \mbox{on $\Omega$}, \s U=0 \s \mbox{on $\partial \Omega$} \,.
\eeq
Multiply this by $\xi _1 \p _1-U $, and integrate
\beq
\lbl{22abc}
\int_\Omega |\nabla U|^2 \, dx-\la _1 \int_{\Omega} U^2 \,dx+\int_\Omega a(x) \left(\xi _1^2 \p _1^2-U^2 \right) \, dx-\xi _1 \mu _1
\eeq
\[
=-\int_\Omega eU \, dx \,.
\]
Dropping two non-negative terms on the left, we get an estimate from below, leading to
\[
\left(\la _2-\la _1-\gamma \right)\int_\Omega  U^2 \, dx  \leq -\int_\Omega eU \, dx \,.
\]
From this we get an estimate on $\int_\Omega  U^2 \, dx$, and then on $\int_\Omega |\nabla U|^2 \, dx$ from (\ref{22abc}).
\medskip

In case the conditions (\ref{22.2a}) and (\ref{22.2b}) hold, we drop the  terms 
\[
\int_\Omega a(x) \xi _1^2 \p _1^2 \, dx -\xi _1 \mu _1>0
\]
on the left in (\ref{22abc}), and proceed the same way.
\epf

We may assume that $g(0)=0$ in (\ref{n1}), without loss of generality (by expanding $g(0)=a \p _1+\bar e(x)$, and redefining $\mu _1$ and $e(x)$).

\begin{thm}
Assume that $e(x) \in L^2(\Omega)$, $g(u) \in C^1(R)$ and $g(0)=0$, $g'(u)<\la _2$ for all  $u \in (-\infty,\infty)$,  and moreover there exist constants $\gamma _1,\gamma _2$, with $\la _1 <\gamma _1<\gamma _2<\la _2$ so that 
\beq
\lbl{g8}
\la _1 <\gamma _1<\frac{g(u)}{u}<\gamma _2<\la _2 \,, \s\s \mbox{for all  $u \in (-\infty,\infty)$} \,.
\eeq
Then all solutions of (\ref{n1}) lie on a unique continuous solution curve $(u(x),\mu _1)(\xi _1)$, and  $\lim _{\xi _1 \ra -\infty} \mu _1(\xi _1)=-\infty$,
$\lim _{\xi _1 \ra \infty} \mu _1(\xi _1)=\infty$, so that the problem (\ref{n1}) is solvable for any $f(x) \in L^2(\Omega)$. If  $g'(u)$ is either strictly increasing or strictly decreasing for all $u \in R$,  then the problem (\ref{n1}) has a unique solution.
\end{thm}

\pf
Observe that (\ref{g8}) implies that the conditions (\ref{n2}) and (\ref{n2a}) hold, and hence by Theorem \ref{thm:1} there exists
 a unique  continuous solution curve $(u(x),\mu _1)(\xi _1)$ for the problem (\ref{n1}). We claim that $\mu _1(\xi _1)$ cannot remain bounded from above as $\xi _1 \ra \infty$, and hence $\lim _{\xi _1 \ra \infty} \mu _1(\xi _1)=\infty$. Indeed, in such a case the condition (\ref{22.2b}) holds, and then by Lemma \ref{lma:6} $||U||_{L^2(\Omega)}$ is bounded uniformly in $\xi _1$. By the elliptic regularity $||U||_{C^1(\bar \Omega)}$ is bounded uniformly in $\xi _1$, and hence $u(x)>0$ for large $\xi _1>0$. Then as $\xi _1 \ra \infty$
\[
\mu _1=\int _{\Omega} g(u) \p _1 \, dx-\la _1 \xi _1 \geq \gamma _1 \int _{\Omega} u \p _1 \, dx-\la _1 \xi _1=\left( \gamma _1-\la _1 \right)\xi _1 \ra \infty \,,
\]
a contradiction. Similarly one shows that $\lim _{\xi _1 \ra -\infty} \mu _1(\xi _1)=-\infty$.
\medskip

In case $g'(u)$ is strictly monotone, the solution curve $\mu _1(\xi _1)$ cannot have any turns, because a turn  would result in at least three solutions of (\ref{n1}) at some value of $\mu _1$, contradicting Theorem \ref{thm:g3}.
\epf

The following result (which was already stated in A. Ambrosetti and G.  Prodi \cite{A}, p.$163$) shows another  use of the principal solution curve. Observe that  no assumption is made on the order of subsolution  and supersolution.
\begin{prop}
Assume that the condition (\ref{n2}) holds, and moreover $g(u)=\la u+h(u)$, with $0 \leq \la <\la _2$, and $\lim _{|u| \ra \infty} \frac{h(u)}{u}=0$. Assume also that the problem 
\beq
\lbl{n3}
\Delta u +g(u)=0 \s \mbox{on $\Omega$}, \s u=0 \s \mbox{on $\partial \Omega$} \,,
\eeq
has a subsolution $\psi (x)$ and supersolution $\p (x)$ (without requiring that $\psi \leq \phi$). Then the problem (\ref{n3}) has a solution.
\end{prop}

\pf
Embed   the problem (\ref{n3}) into
\beq
\lbl{n4}
\Delta u +g(u)=\mu _1 \p _1 \s \mbox{on $\Omega$}, \s u=0 \s \mbox{on $\partial \Omega$} \,,
\eeq
and consider  the continuous solution curve $(u(x),\mu _1)(\xi _1)$, given by Theorem \ref{thm:1}. Suffices to show that $\mu _1(\xi _1)$ changes sign (then at $\mu _1(\xi^0 _1)=0$ we obtain a solution of (\ref{n3})). Assume, on the contrary, that say $\mu_1(\xi _1)>0$ for all $\xi _1$. Then solutions of (\ref{n4}) are subsolutions of (\ref{n3}). Since $u(x)=\xi _1 \p _1(x)+U(x)$, with $U(x)/\xi _1$ uniformly small by the Proposition \ref{prop:inf} and the elliptic regularity, we have a continuous family of subsolutions of  (\ref{n3}) extending from $-\infty$ to $\infty$, uniformly in $x \in \Omega$. But then it is impossible for a supersolution $\p (x)$ to exist. (By the strong maximum principle we obtain a contradiction at a point where a subsolution touches from below the supersolution $\p (x)$. This type of argument is sometimes referred to as Serrin's sweeping principle.) Alternatively, the proof could be completed by observing that we can produce an ordered pair of a subsolution and a supersolution.
\epf

In \cite{K0} we studied a case where a supersolution was below a subsolution, and it was possible to construct a sequence of monotone increasing iterations beginning with a supersolution, and  monotone decreasing iterations beginning with a subsolution. This result was  generalized by J. Shi \cite{S}.

\medskip

We now obtain multiplicity results for the problem  (\ref{n3}), by embedding it into  (\ref{n4}).  

\begin{prop}\lbl{propp}
For the  problem (\ref{n3}) assume that the condition (\ref{n2}) holds, and
\beq
\lbl{n4a}
g(0)=0 \,, \s g'(0)<\la _1 \,,
\eeq
\beq
\lbl{n4b}
\frac{g(u)}{u} \geq \gamma>\la _1 \,, \s \s \mbox{for $|u|>\rho$} \,,
\eeq
for some constants $\gamma$ and $\rho$. Then the  problem (\ref{n3}) has a positive solution and a negative solution. If, moreover, the function $\frac{g(u)}{u}$ is decreasing on $(-\infty,0)$ and increasing on $(0,\infty)$ then  the  problem (\ref{n3}) has exactly two solutions. All solutions of  the  problem (\ref{n3}) can be numerically computed by following  the solution curve $(u(x),\mu _1)(\xi _1)$ of (\ref{n4}) starting with the trivial solution at $\xi _1=0$. 
\end{prop}

\pf
We consider again the continuous solution curve $(u(x),\mu _1)(\xi _1)$ of (\ref{n4}), given by Theorem \ref{thm:1} (see Remark \ref{rmk:1}). Clearly, $\mu _1 (0)=0$, corresponding to the trivial solution. From  (\ref{n4})
\beq
\lbl{n5}
\mu _1(\xi _1)=-\la _1 \xi _1+\int _{\Omega} g(u(x)) \p _1(x) \, dx \,.
\eeq
We claim that $\mu _1(\xi _1)<0$ ($>0$) for $\xi _1>0$ ($<0$) and small.  For $u$ small, $g(u)=g'(0)u+O(u^2)$, and then
\beq
\lbl{n6}
\mu _1(\xi _1)=\left(g'(0)-\la _1 \right)\xi _1+O \left(\int _{\Omega} u^2 \p _1 \, dx \right)
\eeq
The claim will follow once we show that $\int _{\Omega} u^2 \p _1 \, dx=o(\xi _1)$, as $\xi _1 \ra 0$. We have $|g(u)| \leq c_1 |u|$ for some $c_1<\la _1$, on some interval $(-\ep, \ep)$ around $0$. Substituting $u=\xi _1 \p _1+U$ into (\ref{n4}) gives
\[
-\la _1 \xi _1 \p _1+\Delta U+g \left(\xi _1 \p _1+U \right)=\mu _1 \p _1 \,.
\]
Using the continuity of $U(\xi _1)$ and $U(0)=0$, we see that $-\ep<\xi _1 \p _1+U<\ep$ for $|\xi _1|$ small. Then
multiplying by $U$ and integrating
\[
\la _2 \int _{\Omega} U^2 \, dx \leq \int _{\Omega} |\nabla U|^2 \, dx \leq c_1 \int _{\Omega} |\xi _1| \, \p _1  | U| \, dx+c_1 \int _{\Omega} U^2 \, dx \,,
\]
from which one concludes that $\int _{\Omega} U^2 \, dx \leq c_2 \xi _1^2$. Then
\[
\int _{\Omega} u^2 \p _1 \, dx \leq c_3 \int _{\Omega} u^2  \, dx=c_3 \left( \int _{\Omega} U^2 \, dx+\xi _1^2 \right) \leq c_4 \xi _1^2 \,,
\]
and the claim follows by (\ref{n6}).
\medskip

We just saw  that for $\xi _1 >0$ and small, $u(x,\xi _1)$ is small and $\mu _1<0$. Since $g'(0)<\la _1$, it follows by the maximum principle applied to (\ref{n4}) that $u(x,\xi _1)>0$ for $\xi _1 >0$ and small. We claim that $u(x,\xi _1)>0$ so long as $\mu _1<0$. Indeed, the solution $u(x,\xi _1)$ of (\ref{n4}) gives a supersolution of (\ref{n3}), while zero is a solution of (\ref{n3}), and the claim follows by the strong maximum principle. From  (\ref{n4b}), $g(u)>\gamma u-A$ for some $A>0$ when $u>0$, and since $u(x,\xi _1)>0$, it follows by (\ref{n5}) that $\mu _1(\xi _1)>(\gamma-\la _1) \xi _1-A \int _{\Omega} \p _1 \, dx>0$ for large  $\xi _1 >0$, so that a root of $\mu _1(\xi _1)$ must be reached, giving a positive solution of (\ref{n3}). Existence of a  negative  solution of (\ref{n3}) is proved similarly.
\medskip

If $\frac{g(u)}{u}$ is monotone but there are two positive (negative) solutions $u$ and $v$, then from the corresponding equations $\int _{\Omega} uv \left(\frac{g(u)}{u}-\frac{g(v)}{v} \right) \, dx=0$, a contradiction since $u$ and $v$ are ordered by Proposition \ref{prop:g2}.
\epf

\noindent
{\bf Remark} $\s$
One can use a more traditional curve following by embedding the problem (\ref{n3}) into 
\[
\Delta u +\la g(u)=0 \s \mbox{on $\Omega$}, \s u=0 \s \mbox{on $\partial \Omega$} \,.
\]
Let us assume that $\lim _{u \ra \infty} \frac{g(u)}{u}=\gamma $ exists, $\gamma \in (\la _1, \la _2)$. Then there exists a solution curve bifurcating from zero at $\la =\frac{\la _1}{g'(0)}>1$ and going to infinity at $\la =\frac{\la _1}{\gamma}<1$, which passes through $\la =1$, giving a solution of (\ref{n3}). However, this curve may make multiple turns, while there are no turns if continuation in $\xi _1$ is used.
\medskip 

If the conditions at zero and infinity are reversed for the  problem (\ref{n3}), the situation is similar, moreover the existence  of solutions can be proved without the condition (\ref{n2}).

\begin{prop}
For the  problem (\ref{n3}) assume that 
\beq
\lbl{n4c}
g(0)=0 \,, \s g'(0)>\la _1 \,,
\eeq
\beq
\lbl{n4d}
\frac{g(u)}{u} \leq \gamma < \la _1 \,, \s \s \mbox{for $|u|>\rho$} \,,
\eeq
for some constants $\gamma >0$ and $\rho >0$. Then the  problem (\ref{n3}) has a positive solution and a negative solution. 
\medskip

Assume additionally that the condition (\ref{n2}) holds.  Then all solutions of  the  problem (\ref{n3}) can be numerically computed by following  the solution curve $(u(x),\mu _1)(\xi _1)$ of (\ref{n4}) starting with the trivial solution at $\xi _1=0$. If the function $\frac{g(u)}{u}$ is increasing on $(-\infty,0)$ and decreasing on $(0,\infty)$ then  the  problem (\ref{n3}) has exactly two solutions.
\end{prop}

\pf
The function $\epsilon \p _1(x)$ is a subsolution of the  problem (\ref{n3}), for small enough $\epsilon >0$. The condition (\ref{n4d}) implies that $g(u) <\gamma u+c$ for all $u>0$ and $c \geq 0$. Let $\hat \Omega \supset \Omega$ be a slightly larger domain, with the principal eigenpair $(\hat \la _1, \hat \p _1(x))$ such that $\hat \la < \gamma$. Then $\hat \p _1(x))>0$ for $x \in \bar \Omega$, and the function $M \hat \p _1(x)$ is a supersolution of the  problem (\ref{n3}), for large  enough $M$, proving the existence of positive solution. Similarly we use   a subsolution $-M \hat \p _1(x))$ and a supersolution $-\epsilon \p _1(x)$, with small $\epsilon >0$ and large $M>0$, to prove the existence of a negative solution.
\medskip

If  the condition (\ref{n2}) holds, similarly to Proposition \ref{propp} we have a solution curve $(u(x),\mu _1)(\xi _1)$ of (\ref{n4}) with $\mu _1(\xi _1)>0$ ($<0$) for $\xi _1>0$ ($<0$) and small, and $\mu _1(\xi _1)$ changing sign for $|\xi _1|$ large.
\epf

\section{Two important classes of  equations}
\setcounter{equation}{0}
\setcounter{thm}{0}
\setcounter{lma}{0}
\setcounter{rmk}{0}
\setcounter{prop}{0}
\setcounter{cor}{0}

We now give a detailed result for the problem (\ref{80a}) below. It turns out that the case of resonance, when $\la =\la _1$, is not particularly distinguished from the other $\la$'s.
\begin{thm}\lbl{thm5.1!}
Consider the problem
\beq
\lbl{80a}
\Delta u+\la u-u^3=f(x) \,, \, \s \mbox{for $x \in \Omega$} \,,\, \s u=0 \s \mbox{on $\partial \Omega$} \,,
\eeq
with a parameter $\la \in (0,\la _2)$, and $f(x) \in L^2(\Omega)$. Decompose $f(x)=\mu _1 \p _1+e(x)$, with  $\int _{\Omega} e(x) \p _1(x) \, dx=0$, and  $u(x)=\xi _1 \p _1(x)+U(x)$, with $\int _{\Omega} U(x) \p _1(x) \, dx=0$, and assume that $|e(x)| \leq M$ for some $M>0$ and $x \in \Omega$. Then for each $\xi _1\in R$ one can find a unique solution pair $(u(x),\mu _1)$, and the solution curve $(u(x),\mu _1)(\xi _1)$ is continuous. Moreover,  as $\xi _1 \ra \pm \infty$, we have $\mu _1 \ra \mp \infty$ and $\frac{u(x)}{\xi _1} \ra \p _1(x)$ in $H^1(\Omega)$. The curve $\mu _1=\mu _1 (\xi _1)$ is decreasing for $|\xi _1|$  (or for $|\mu _1|$) large (implying the uniqueness of  solution of (\ref{80a})), and the solution $u(x)$ is positive (negative) for $\xi _1 >0$ ($\xi _1 <0$) and large. In case $\la \in (0,\la _1]$, the curve $\mu _1=\mu _1 (\xi _1)$ is decreasing for all $\xi _1$. In case $\la \in (\la _1,\la _2)$, the curve $\mu _1=\mu _1 (\xi _1)$ makes at least two turns for $e(x)$ small.
\end{thm}

\pf
Proposition \ref{prop:1} provides us with the unique solution curve, we only need to prove its properties. Let us begin with curve's behavior as $\xi _1 \ra \pm \infty$. We claim that $|\mu _1| \ra \infty$ as $\xi _1 \ra \pm \infty$. Indeed, if we assume that $|\mu _1|$ is bounded, then $f(x)=\mu _1 \p _1+e(x)$ is bounded in $L^2(\Omega)$. Multiplying (\ref{80a}) by $u$ and integrating, obtain
\beq
\lbl{80b}
\int_\Omega |\nabla u|^2 \, dx+\int_\Omega u^4 \, dx-\la \int_\Omega u^2 \, dx=-\int_\Omega uf(x) \, dx \,.
\eeq
We conclude the boundness of $\int_\Omega u^2 \, dx$ (the $\int_\Omega u^4 \, dx$ term controls $\int_\Omega u^2 \, dx$). Using 
\beq
\lbl{80c}
\int_\Omega u^2 \, dx=\xi _1 ^2\int_\Omega \p_1 ^2 \,dx+\int_\Omega U^2 \, dx \,, 
\eeq
we conclude that $|\xi _1|$ is bounded, a contradiction. Assume, contrary to what we wish to prove, that $\mu _1 \ra \pm \infty$,  as $\xi _1 \ra \pm \infty$. Write (\ref{80b}) in the form
\[
\int_\Omega |\nabla u|^2 \, dx+\int_\Omega u^4 \, dx-\la \int_\Omega u^2 \, dx=-\xi_1 \mu _1 - \int_\Omega U e \, dx \,.
\]
Dropping the negative term $\xi_1 \mu _1$, and estimating $|\int_\Omega U e \, dx|  \leq \ep \int_\Omega u^2 \, dx+c(\ep) \int_\Omega  e^2 \, dx$, we obtain again a bound from above on $\int_\Omega u^2 \, dx$, which implies that $|\xi _1|$ is bounded by (\ref{80c}), a contradiction. By Proposition \ref{prop:1} the estimate (\ref{18}) holds, so that  $\frac{u(x)}{\xi _1}= \p _1(x)+\frac{U(x)}{\xi _1}\ra \p _1(x)$ in $H^1(\Omega)$ as $\xi _1 \ra \pm \infty$.
\medskip

We claim next that the solutions of (\ref{80a}) are positive for $\xi _1 >0$ and large, or when $\mu _1<0$ and large in absolute value. Denote $\Omega _-= \{ x \in \Omega \, |  \, u(x)<0 \}$. Then from  (\ref{80a}), when $|\mu _1|$ is large enough,  
\[
\Delta u=-\la u+u^3+\mu _1 \p _1 (x)+e(x)<0 \s \mbox{on $\Omega _-$} \,, u=0 \s \mbox{on $\partial \Omega _-$} \,,
\]
which implies that $u(x)>0$ on $\Omega _-$, a contradiction. (It is only at this step that we use the boundness of $e(x)$.) Since $u(x)>0$, the nonlinear term in  (\ref{80a}) is concave for  $\xi _1 >0$ and large. It follows by the Theorem \ref{thm:g3} that the problem  (\ref{80a})  has at most two solutions. But then the solution curve cannot turn, and has to remain decreasing in $\xi _1$. Indeed, if the curve turned, it would have to turn again (since $\mu _1 \ra - \infty$ as $\xi _1 \ra \infty$), giving us three solutions of  (\ref{80a}), a contradiction.
\medskip

We show next that for $\la \leq \la _1$ the solution of (\ref{80a}) is  unique, which implies that the solution curve $\mu _1=\mu _1 (\xi _1)$ is decreasing for all $\xi _1 \in R$. Indeed, if $u(x)$ and $v(x)$ are two solutions of  (\ref{80a}), then $w=u-v$ satisfies
\beq
\lbl{83}
\Delta w+\la w  -\left(u^2+uv+v^2 \right)w=0 \,,  \s x \in \Omega, \s w=0 \;\; \mbox{on $\partial \Omega$} \,.
\eeq
Since $u^2+uv+v^2$ is non-negative and non-zero, we conclude that $w \equiv 0$.
\medskip

We now turn to the  case when $\la>\la _1$ and $e(x)$ is small. Start by assuming that  $e(x)=0$.  Then $\mu _1(0)=0$, and $\mu _1(\xi _1)>0$ for $\xi _1>0$ and small, as seen by multiplying  (\ref{80a}) by $\p _1$ and integrating. Since the function $\mu _1(\xi _1)$ is eventually decreasing, it must have a local maximum for some $\xi _1>0$. Similarly, there is a local minimum  for some $\xi _1<0$. The solution curve remains similar for small $e(x)$. 
\epf

The problem (\ref{80a}) was studied in a series of papers by P.T. Church, J.G. Timourian and their coworkers, see e.g., \cite{B}, \cite{D} and the references therein.
The result of P.T. Church et al \cite{D} is more detailed in many (but not all) respects than our Theorem \ref{thm5.1!}, but it only applies for $n \leq 4$. The same restriction on the dimension appeared also in M.S. Berger et al \cite{B}, and other papers.
\medskip

The $u^3$ term in (\ref{80a}) can be changed to $u |u|^{p-1}$ with $p>1$, although some arguments (like the one in (\ref{83})) would need to be modified. However, we cannot replace $u^3$ term by $u^2$, to handle the logistic equation with harvesting, for which a similar result was established in Y. Wang et al \cite{W}.
\medskip

We now consider a class resonant  problem
\beq
\lbl{23}
\Delta u +\la _1 u+g(u)=\mu _1 \p _1+e(x) \s \mbox{on $\Omega$}, \s u=0 \s \mbox{on $\partial \Omega$} \,,
\eeq
with $e(x) \in \p _1 ^\perp$ in $L^2(\Omega)$. We wish to find a solution pair $(u, \mu _1)$. The famous  existence result of  E.M. Landesman and A.C. Lazer \cite{L} required $g(u)$ to have finite limits at $\pm \infty$. Then  for 
{\em bounded} $g(u)$, satisfying $ug(u) \geq 0$ for all $u \in R$, and  $\mu _1=0$,
 D.G. de Figueiredo and W.-M.  Ni \cite{FN} have proved the existence of solutions. R. Iannacci, M.N. Nkashama and J.R. Ward \cite{I} generalized this result
to unbounded $g(u)$ satisfying $g'(u) \leq \gamma <\la _2-\la _1$, while still assuming  $\mu _1=0$. An overview of these results can be found in P. Korman \cite{Knew}.
In \cite{K9} we extended  the result of R. Iannacci et al \cite{I} to the problem (\ref{23}), with  $\mu _1 \ne 0$.
We now present a generalization of our result in \cite{K9}, dropping the technical condition (3.2) of that paper.
\begin{thm}\lbl{thm4}
Assume that $g(u) \in C^1(R)$ satisfies
\beq
\lbl{24}
u g(u) >0 \s\s \mbox{for all $u \in R$} \,,
\eeq
\beq
\lbl{25}
 g'(u) \leq \gamma< \la _2-\la _1 \s\s \mbox{for all $u \in R$} \,.
\eeq
Then there is a continuous curve of solutions of (\ref{23}): $(u(\xi _1),\mu _1(\xi _1))$, $u \in H^2(\Omega) \cap H^1_0(\Omega)$, with $-\infty<\xi _1<\infty$, and $\int _\Omega u(\xi _1) \p _1 \, dx=\xi _1$. This curve exhausts the solution set of (\ref{23}). The continuous function $\mu _1(\xi _1)$ is positive for $\xi _1 >0$ and large, and $ \mu _1(\xi _1)<0$ for $\xi _1 <0$ and $|\xi _1|$ large. In particular, $\mu _1(\xi^0 _1)=0$ at some $\xi^0 _1$, concluding the existence of  solution for 
\[
\Delta u +\la _1 u+g(u)=e(x) \s \mbox{on $\Omega$}, \s u=0 \s \mbox{on $\partial \Omega$} \,.
\]
\end{thm}

\pf
By the Theorems \ref{thm:1} and \ref{thm:2} there exists a curve of solutions of (\ref{23}) $(u(\xi _1),\mu _1(\xi _1))$, which 
exhausts the solution set of (\ref{23}). The properties of this curve follow the same way as in \cite{K9}.
\epf

\section{Two classes of oscillatory equations at resonance}
\setcounter{equation}{0}
\setcounter{thm}{0}
\setcounter{lma}{0}
\setcounter{rmk}{0}
\setcounter{prop}{0}
\setcounter{cor}{0}

We now use stationary phase method to obtain  more detailed results in the one dimensional case. In particular, we make use of the $k$-th solution curves.  Recall  
 the following asymptotic formula, see Y.V. Sidorov et al \cite{SFS}.
\begin{lma}\lbl{lma:o1}
Assume that  $f(x)$ and $g(x)$ are of class $C^2[a,b]$ and $g(x)$ has a unique critical point $x_0$ on $[a,b]$, and moreover $x_0 \in (a,b)$ and $g''(x_0) \ne 0$ (so that $x_0$ gives a global max or global min on $[a,b]$). Then as $\la \ra \infty $ the following asymptotic formula holds
\[
\int _a^b f(x) e^{i \la g(x)} \, dx=e^{i \left[ \la g(x_0) \pm \frac{\pi}{4} \right]} \sqrt{ \frac{2 \pi}{\la |g''(x_0)|}} 
\, f(x_0)+O \left(\frac{1}{\la} \right) \,,
\]
where one takes ``plus" if $g''(x_0)>0$ and ``minus"  if $g''(x_0)<0$.
\end{lma}

We now present a class of Dirichlet problems with infinitely many solutions for any right hand side:
\beqa
\lbl{o1}
& u''+\frac{\pi ^2}{L^2} \, u+h(u) \sin u=\mu _1 \sin \frac{\pi}{L} x+e(x) \,, \s x \in (0,L) \,, \\ \nonumber
& u(0)=u(L)=0 \,. \nonumber
\eeqa
Here $\frac{\pi ^2}{L^2}=\la _1$, the principal eigenvalue of $u''$ on $ (0,L)$ corresponding to $\p _1(x)=\sin \frac{\pi}{L} x$, $\mu _1 \in R$, $e(x) \in C (0,L)$ satisfies $\int _0^L e(x) \sin \frac{\pi}{L} x \, dx=0$. Decompose $u(x)=\xi _1  \sin \frac{\pi}{L} x+U(x)$, with $\xi _1=\frac{2}{L} \int _0^L u(x) \sin \frac{\pi}{L} x \, dx$ and  $\int _0^L U(x) \sin \frac{\pi}{L} x \, dx=0$. (The normalization of $\p _1(x)$ is different from the previous sections, so that $\xi _1$ and $\mu _1$ are changed accordingly by a factor.)

\begin{thm}\lbl{thm:o1}
Assume that $e(x) \in C (0,L)$, and $h(u) \in C^2(R)$  satisfies
\beq
\lbl{o2--}
|h(u)| <\frac{3 \pi ^2}{L^2} |u|+c \,, \s \mbox{for all $u \in R$ and some $c \geq 0$}\,,
\eeq
\beq
\lbl{o2-}
\lim _{u \ra \infty} \frac{h(u)}{u^p}=h_0 \,, \s \mbox{with constants $p \in (\frac{1}{2},1)$ and $h_0>0$} \,.
\eeq
Then for any $\mu _1 \in R$ the problem (\ref{o1}) has infinitely many classical solutions. Moreover, as $\xi _1 \ra \pm \infty$, we have $\frac{u(x)}{\xi _1} \ra \sin \frac{\pi}{L} x$ in $C^2(0,L)$, and 
\beq
\lbl{o2}
\mu _1(\xi _1) \sim \frac{2 \sqrt{2}}{\sqrt{\pi \xi _1}} \sin \left( \xi _1-\frac{\pi}{4} \right)  h\left( \xi _1(1+o(1) \right), \s \mbox{as $\xi _1 \ra  \infty$} \,,
\eeq
\beq
\lbl{o2++}
\mu _1(\xi _1) \sim \frac{2 \sqrt{2}}{\sqrt{\pi |\xi _1|}} \sin \left( \xi _1+\frac{\pi}{4} \right)  h\left( \xi _1(1+o(1) \right), \s \mbox{as $\xi _1 \ra  -\infty$} \,.
\eeq
\end{thm}

\pf
Assume first that $\xi _1>0$.
By (\ref{o2--}), the Proposition \ref{prop:n1} applies, and the problem (\ref{o1}) has  a curve of solutions $(u(\xi _1),\mu _1(\xi _1))$. By the Proposition \ref{prop:inf}, $\frac{u(x)}{\xi _1} \ra \sin \frac{\pi}{L} x$ in $L^2(0,L)$, and by the elliptic regularity $\frac{u(x)}{\xi _1} \ra \sin \frac{\pi}{L} x$ in $C^2(0,L)$, or $u(x)=\xi _1 \sin \frac{\pi}{L} x+\xi _1v(x)$, with $||v(x)||_{C^2(0,L)}=o(1)$, as $\xi _1 \ra \infty$. It follows that the function $u(x)$ is unimodular with the point of maximum lying near $L/2$, and its maximum value is equal to $\xi _1(1+o(1))$, for large $\xi _1$. Multiply the equation (\ref{o1}) by $\sin \frac{\pi}{L} x$, then  integrate over $(0,L)$ and use Lemma \ref{lma:o1}, with $g(x)=u(x)$. As $\xi _1 \ra \infty$, obtain
\[
\mu _1 \, \frac{L}{2}=\int_0^L h(u(x)) \sin u(x) \sin \frac{\pi}{L} x \,dx={\rm Im} \, \int_0^L h(u(x))  \sin \frac{\pi}{L} x \, e^{i u(x)}\,dx
\]
\[
\sim {\rm Im} \, e^{i \left[\xi _1(1+o(1))-\frac{\pi}{4}  \right]}  \sqrt{\frac{2\pi}{\pi^2/L^2 \left[ \xi _1(1+o(1)) \right]} } \, h\left( \xi _1(1+o(1) \right)  \,,
\]
which implies (\ref{o2}).    By our assumptions on $h(u)$, the formula (\ref{o2}) implies that  $\mu _1(\xi _1)$ is oscillatory with the amplitude tending to infinity as $\xi _1 \ra \infty$, and hence it assumes each value in $R$ infinitely many times. In case $\xi _1<0$, Lemma \ref{lma:o1}, with $g(x)=-u(x)$, leads to the asymptotic formula (\ref{o2++}), and again  $\mu _1(\xi _1)$ is oscillatory with the amplitude tending to infinity as $\xi _1 \ra -\infty$. The proof is not finished yet, because the Proposition \ref{prop:n1} did not provide us with the continuity of $\mu _1(\xi _1)$. Fix a value of $\mu _1=\mu _1^0$. Then any solution of (\ref{o1}) with $\mu _1>\mu _1^0$ ($\mu _1<\mu _1^0$) provides us with a subsolution (supersolution) of  (\ref{o1})  at $\mu _1=\mu _1^0$. Since $u(x) \sim \xi _1 \sin \frac{\pi}{L} x$, we can arrange for an ordered subsolution - supersolution pair, providing a solution of (\ref{o1}) between them. (There are infinitely many choices of a supersolution. Select one with a larger $\xi _1$ to exceed any subsolution). After producing one such ordered pair, one starts working on the next one. We conclude the existence of infinitely many solutions of (\ref{o1})  at $\mu _1=\mu _1^0$. 
\epf

\begin{rmk}
\begin{enumerate}
  \item It follows that the problem
\[
u''+\frac{\pi ^2}{L^2} \, u+h(u) \sin u=f(x) \,, \s x \in (0,L) \,, \s 
 u(0)=u(L)=0 
\]
has infinitely many solutions for any continuous $f(x)$.
  \item Since $u(x) \sim \xi _1 \sin \frac{\pi}{L} x$, it follows that the maximum of $u(x)$ is asymptotic to the first harmonic $\xi _1$, for large $\xi _1$.
  \item 
In case $\mu _1=0$, the existence of infinitely many solutions also follows from R. Schaaf  and K. Schmitt \cite{SS}, see also D. Costa et al \cite{CJSS}.
\end{enumerate}
\end{rmk}

\noindent
{\bf Example} $\;$ We computed the solution curve $\mu _1= \mu _1(\xi _1)$ for the following example
\beqa 
\lbl{o2a}
& \s\s u''+\pi ^2\, u+5 \left(u^2+1 \right)^{\frac{5}{12}} \sin u=\mu _1 \sin \pi x+0.2 \sin 2\pi x \,, \s x \in (0,1) \,, \\ \nonumber
& u(0)=u(1)=0 \,. \nonumber
\eeqa

\noindent
Here $\la _1=\pi ^2$, $\p _1(x)=\sin \pi x$, and $e(x)=0.2 \sin 2\pi x \perp \p _1(x)$ in $L^2(0,1)$.
The solution curve $\mu _1= \mu _1(\xi _1)$ (solid line) is presented in Figure \ref{fig:1}. Notice an excellent agreement with the asymptotic formula (\ref{o2}) (dashed line).

\begin{figure}
\begin{center}
\scalebox{0.95}{\includegraphics{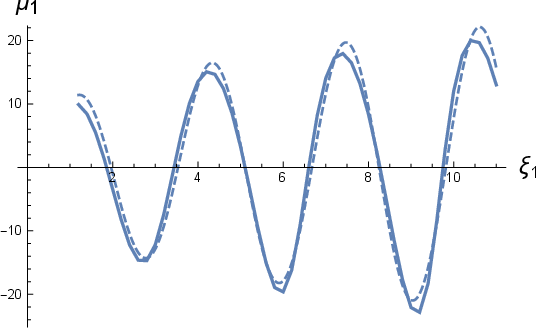}}
\end{center}
\caption{ The  solution curve $\mu _1= \mu _1(\xi _1)$ of (\ref{o2a}), compared with (\ref{o2})}
\lbl{fig:1}
\end{figure}
\medskip

\begin{thm}
Assume that (\ref{o2-}) holds with $p \in (0,\frac{1}{2})$. Then the problem (\ref{o1}) has infinitely many solutions for $\mu _1=0$, and at most finitely many solutions for $\mu _1 \ne 0$, on any solution curve. As $\xi _1 \ra \pm \infty$, the asymptotic formulas (\ref{o2}) and (\ref{o2++}) hold, and moreover
\beq
\lbl{o3a}
u(x) \sim \xi _1 \sin \frac{\pi}{L} x+E(x) \,,
\eeq
where $E(x)$ is the unique solution of 
\beq
\lbl{o3b}
u''+\frac{\pi ^2}{L^2} \, u=e(x) \,, \s    u(0)=u(L)=0 \,, \int_0^L u(x) \sin \frac{\pi}{L} x \, dx=0 \,.
\eeq
\end{thm}

\pf
As above we have a solution curve $(u(x), \mu _1)(\xi _1)$, and  the asymptotic formulas (\ref{o2}) and (\ref{o2++}) hold, which implies that $\mu _1(\xi _1) \ra 0$ as $\xi \ra \pm \infty$, justifying the multiplicity claims. Let $G(x,z)$ be the (bounded) modified Green's function for (\ref{o3b}), see e.g.,  p. 136 in E.C. Young \cite{Y}. Express the solution of  (\ref{o1}) as
\[
u(x)=-\int_0^L G(x,z) h(u(z)) \sin u(z) \, dz+\mu _1 \int_0^L G(x,z) \sin \frac{\pi}{L} z \, dz
\]
\[
+E(x)+c\sin \frac{\pi}{L} x \,.
\]
As $\xi _1 \ra \pm \infty$, we have $\mu _1(\xi _1) \ra 0$, and the first integral on the right tends to zero by a similar argument, based on Lemma \ref{lma:o1}. It follows that $u(x) \ra E(x)+c\sin \frac{\pi}{L} x$, and then $c=\xi _1$, since $\frac{u(x)}{\xi _1} \ra \sin \frac{\pi}{L} x$ as $\xi _1 \ra \pm \infty$.
\epf

For the case of higher eigenvalues we restrict to a model problem ($k>0$ is an integer)
\beqa
\lbl{o4}
& u''+\frac{k^2 \pi ^2}{L^2} \, u+ \sin u=\mu _k \sin \frac{k\pi}{L} x+e(x) \,, \s x \in (0,L) \,, \\ \nonumber
& u(0)=u(L)=0 \,, \nonumber
\eeqa
although our results can be easily generalized in a number of ways. 
Here $\frac{k^2\pi ^2}{L^2}=\la _k$, the $k$-th eigenvalue of $u''$ on $ (0,L)$ with zero boundary conditions, $\mu _k \in R$, $e(x) \in C (0,L)$ satisfies $\int _0^L e(x) \sin \frac{k\pi}{L} x \, dx=0$. Decompose $u(x)=\xi _k  \sin \frac{k\pi}{L} x+U(x)$, with $\int _0^L U(x) \sin \frac{k\pi}{L} x \, dx=0$, as  above.

\begin{thm}\lbl{thm:o2}
Assume that the constants $k \in N$ and $L \in R$ satisfy
\beq
\lbl{o5}
\frac{(k-1)^2 \pi ^2}{L^2}+1 <\frac{k^2 \pi ^2}{L^2}<\frac{(k+1)^2 \pi ^2}{L^2}-1 \,.
\eeq
Then  the problem (\ref{o4}) has a unique continuous solution curve $(u(x),\mu _k)(\xi _k)$, for all $\xi _k \in R$. Moreover, as $\xi _k \ra \pm \infty$, we have $\frac{u(x)}{\xi _k} \ra \sin \frac{k\pi}{L} x$ in $C^2(0,L)$, and 
\beq
\lbl{o6}
 \mu _k(\xi _k) \sim 2 \sqrt{ \frac{2}{\pi \xi _k}} \sin \left(  \xi _k-\frac{\pi}{4} \right), \mbox{as $\xi _k \ra  \infty$} \,,
\eeq
\beq
\lbl{o6a}
 \mu _k(\xi _k) \sim 2 \sqrt{ \frac{2}{\pi |\xi _k|}} \sin \left(  \xi _k+\frac{\pi}{4} \right), \mbox{as $\xi _k \ra  -\infty$} \,.
\eeq

It follows that when $\mu _k=0$ the problem (\ref{o4}) has infinitely many solutions, there are at most finitely many solutions for any $\mu _k \ne 0$, and there are no solutions for $|\mu _k|$ large enough.
\end{thm}

\pf
The condition \ref{o5} ensures that the function $g(u)=\frac{k^2 \pi ^2}{L^2} \, u+ \sin u$ satisfies $\la _{k-1}<g'(u)<\la _{k+1}$. By the Remark  \ref{rmk:*} all solutions of (\ref{o4}) lie on a unique continuous solution curve $(u(x),\mu _k)(\xi _k)$. By Proposition \ref{prop:infa}, $\frac{u(x)}{\xi _k} \ra \sin \frac{k\pi}{L} x$ in $H^2(0,L)$, and by the elliptic regularity $\frac{u(x)}{\xi _k} \ra \sin \frac{k\pi}{L} x$ in $C^2(0,L)$, as $\xi _k \ra \pm \infty$. It follows that  the function $u(x)/\xi _k$ has the the same number of points of local maximums and minimums as that of $ \sin \frac{k\pi}{L} x$, and that these points as well as the roots of $u(x)/\xi _k$ tend to the corresponding points of $ \sin \frac{k\pi}{L} x$, as $\xi _k \ra \pm \infty$.
Multiply the equation (\ref{o4}) by $\sin \frac{k\pi}{L} x$, then  integrate over $(0,L)$. As $\xi _k \ra \infty$, similarly to the Theorem \ref{thm:o1} obtain
\beq
\lbl{o7}
\mu _k(\xi _k) \frac{L}{2} \sim {\rm Im} \, \int_0^L   \sin \frac{k \pi}{L} x \, e^{i \xi _k \sin \frac{k\pi}{L} x}\,dx  \,.
\eeq

The case $k=1$ was covered by the formula (\ref{o2}) above (when $h(x) \equiv 1$), so assume that $k \geq 2$. The function $\sin \frac{k \pi}{L} x$ has its first root at $\frac{L}{k}$ and second one at $\frac{2L}{k}$, it is positive on $(0,\frac{L}{k})$ and negative on $(\frac{L}{k},\frac{2L}{k})$. By Lemma \ref{lma:o1}, with $g(x)=\sin \frac{k \pi}{L} x$ and $|g''(x_0)|=\frac{k^2 \pi ^2}{ L^2}$,
\[
 \int_0^{\frac{L}{k}}   \sin \frac{k \pi}{L} x \, e^{i \xi _k \sin \frac{k\pi}{L} x}\,dx \sim \frac{L}{k} \sqrt{\frac{2}{\pi \xi _k} }e^{i \left(\xi _k -\pi/4 \right)} \,.
\]
Similarly, over the negative hump
\[
 \int_{\frac{L}{k}}^{\frac{2L}{k}}   \sin \frac{k \pi}{L} x \, e^{i \xi _k \sin \frac{k\pi}{L} x}\,dx \sim -\frac{L}{k} \sqrt{\frac{2}{\pi \xi _k}}e^{i \left(-\xi _k +\pi/4 \right)} \,.
\]
Adding these, we see that over the first pair of humps
\[
{\rm Im}  \int_0^{\frac{2L}{k}}   \sin \frac{k \pi}{L} x \, e^{i \xi _k \sin \frac{k\pi}{L} x}\,dx \sim \frac{2L}{k} \sqrt{\frac{2}{\pi \xi _k}} \sin \left( \xi _k-\frac{\pi}{4} \right) \,.
\]
If $k$ is even, there are $k/2$ such pairs of humps, and then 
\[
{\rm Im}  \int_0^L   \sin \frac{k \pi}{L} x \, e^{i \xi _k \sin \frac{k\pi}{L} x}\,dx \sim L \sqrt{\frac{2}{\pi \xi _k}} \sin \left( \xi _k-\frac{\pi}{4} \right) \,,
\]
which implies the first formula in (\ref{o6}). In case $k$ is odd, there are $\frac{k-1}{2}$ pairs of humps plus one more positive hump over $(\frac{(k-1)L}{k},L)$. By Lemma \ref{lma:o1}, the last positive hump contributes (the same as for all other positive humps)
\beq
 {\rm Im}  \int_{_\frac{(k-1)L}{k}}^L   \sin \frac{k \pi}{L} x \, e^{i \xi _k \sin \frac{k\pi}{L} x}\,dx \sim \frac{L}{k} \sqrt{ \frac{2}{\pi \xi _k}}  \sin ( \xi _k-\pi/4) \,.
\eeq
Adding to that the contribution from $\frac{k-1}{2}$ pairs of humps, gives
\beqa \nonumber
& {\rm Im}  \int_0^L   \sin \frac{k \pi}{L} x \, e^{i \xi _k \sin \frac{k\pi}{L} x}\,dx \sim \frac{k-1}{2}\frac{2L}{k} \sqrt{\frac{2}{\pi \xi _k}} \sin \left( \xi _k-\frac{\pi}{4} \right) +\frac{L}{k} \sqrt{ \frac{2}{\pi \xi _k}}  \sin ( \xi _k-\pi/4)\\ \nonumber
& =L \sqrt{ \frac{2}{\pi \xi _k}} \sin \left(  \xi _k-\frac{\pi}{4} \right) \,,\nonumber
\eeqa
which again leads to the  formula (\ref{o6}). The formula (\ref{o6a}) is established similarly.
\epf

\noindent
{\bf Example }  Using {\em Mathematica}, we computed the solution curve $\mu _7= \mu _7(\xi _7)$ for the following example of resonance at a higher eigenvalue ($\la _7=49\pi ^2$ on $(0,1)$; the computer program is described in \cite{KS})
\beqa 
\lbl{o7a}
& \s\s\s\s u''+49\pi ^2\, u+ \sin u=\mu _7 \sin 7\pi x+ \sin 3\pi x -2\sin 4\pi x\,, \s x \in (0,1) \,, \\ \nonumber
& u(0)=u(1)=0 \,. \nonumber
\eeqa

\noindent
Observe that here $e(x)=\sin 3\pi x -2\sin 4\pi x \perp \sin 7\pi x$ in $L^2(0,1)$. The solution curve $\mu _7= \mu _7(\xi _7)$ (solid line) is presented in Figure \ref{fig:2}. Notice a remarkable agreement with the asymptotic formula (\ref{o6}) (dashed line, almost indistinguishable). At $\mu _7=0$ the problem (\ref{o7a}) has infinitely many solutions. There are at most finitely many solutions for any $\mu _7 \ne 0$, and there are no solutions for $|\mu _7|$ sufficiently large.

\begin{figure}
\begin{center}
\scalebox{0.95}{\includegraphics{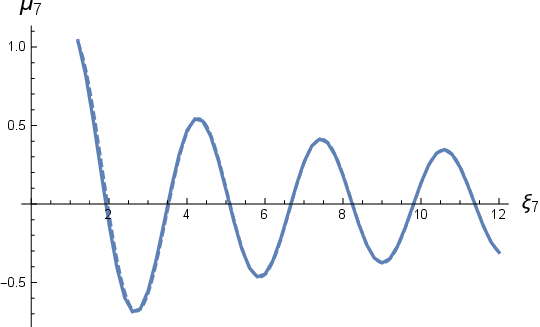}}
\end{center}
\caption{ The  solution curve $\mu _7= \mu _7(\xi _7)$ of (\ref{o7a}), compared with (\ref{o6})}
\lbl{fig:2}
\end{figure}

\end{document}